\documentclass[11pt]{article}
\usepackage{amsmath,amssymb}

\def\r{\rightarrow}
\def\qed{\hfill\vrule height5pt width5pt depth0pt}
\def\one #1{1_{\{#1\}}}

\def\1{\mbox{\bf 1}}

\def\n{\nonumber \\}
\def\phi{\varphi}

\newcommand{\proof}{\noindent {\bf Proof:\ }}

\newtheorem{Theorem}{Theorem}

\begin{document}

\title{The class of distributions associated with the generalized Pollaczek-Khinchine formula}
\author{Offer Kella\thanks{Department of Statistics; The Hebrew University of Jerusalem; Mount Scopus,
Jerusalem 91905; Israel ({\tt offer.kella@huji.ac.il})}
\thanks{Supported in part by grant 434/09 from the Israel Science
Foundation and the Vigevani Chair in Statistics.} }
\date{November 29, 2011}
\maketitle

\abstract{The goal is to identify the class of distributions to which the distribution of the maximum of a L\'evy process with no negative jumps and negative mean (equivalently, the stationary distribution of the reflected process) belongs. An explicit new distributional identity is obtained for the case where the L\'evy process is an independent sum of a Brownian motion and a general subordinator (nondecreasing L\'evy process) in terms of a geometrically distributed sum of independent random variables. This generalizes both the distributional form of the standard Pollaczeck-Khinchine formula for stationary workload distribution in the M/G/1 queue and the exponential stationary distribution of a reflected Brownian motion.}

\vspace{0.1in}
Keywords: Generalized Pollaczek-Khinchine formula, L\'evy process with no negative jumps, spectrally positive L\'evy process, reflected L\'evy process, supremum of a L\'evy process.

\vspace{0.1in}

AMS 2000 Subject Classification: 60G51, 60K25.

\section{Introduction and preliminaries}\label{Sec:intro}
Let $X=\{X_t|\ t\ge0\}$ be a L\'evy process with no negative jumps. It is standard knowledge that in this case $Ee^{-\alpha X_t}=e^{\phi(\alpha)t}$ is finite for all $\alpha\ge 0$ and that
\begin{equation}
\phi(\alpha)=b\alpha+\frac{\sigma^2}{2}\alpha^2+\int_{(0,\infty}\left(e^{-\alpha x}-1+\alpha x\one{x\le 1}\right)\nu(dx)
\end{equation}
where $b$ is real, $\sigma$ is nonnegative and $\nu$ is a measure satisfying
\begin{equation}
\int_{(0,\infty)}(x^2\wedge 1)\nu(dx)<\infty
\end{equation}
where $a\wedge b=\min(a,b)$.

It is also well known that $-EX_t/t=\phi'(0)=-b+\int_{(1,\infty)}x\nu(dx)$ and in particular $EX_t$ is well defined and can be either finite or $+\infty$ but never $-\infty$. In particular if $EX_t<0$ (equivalently $\phi'(0)>0$) then $\int_{(1,\infty)}x\nu(dx)<b$ so that in particular $b>0$ and $\int_{(1,\infty)}x\nu(dx)<\infty$. Therefore, in this case $\phi$ has the following form
\begin{equation}\label{ea:mu}
\phi(\alpha)=\mu\alpha+\frac{\sigma^2}{2}\alpha^2+\int_{(0,\infty}\left(e^{-\alpha x}-1+\alpha x\right)\nu(dx)
\end{equation}
where $\mu=\phi'(0)=b-\int_{(1,\infty)}x\nu(dx)>0$.

When in addition $\int_{(0,1]}x\nu(dx)<\infty$ then the L\'evy process is an independent sum of a Brownian motion and a pure jump subordinator (nondecreasing L\'evy process). In this case, $\phi$ becomes
\begin{equation}
\phi(\alpha)=c\alpha+\frac{\sigma^2}{2}\alpha^2-\int_{(0,\infty)}\left(1-e^{-\alpha x}\right)\nu(dx)\ ,
\end{equation}
where $c=\mu+\int_{(0,\infty)}x\nu(dx)$.

It is also well known that for any L\'evy process with no negative jumps and $\phi'(0)>0$, if we denote $M=\sup_{t\ge0}X_t$, then
\begin{equation}
Ee^{-\alpha M}=\frac{\alpha\phi'(0)}{\phi(\alpha)}
\end{equation}
and that the limiting and stationary distributions of the (Markov) process $W_t=X_t+L_t$, where $L_t=-\inf_{0\le s\le t}X_s$, is the distribution of $M$. We (and others) call the formula $\frac{\alpha\phi'(0)}{\phi(\alpha)}$ the {\em generalized Pollaczek-Khinchine (PK) formula}. The reason for the name is that for the special case of the M/G/1 queue, the underlying driving (L\'evy) process is a compound Poisson process (a subordinator) minus $t$ for which the Laplace-Stieltjes transform of the limiting and stationary distribution of the workload process is the celebrated Pollaczek-Khinchine formula and may be found in virtually all basic queuing theory textbooks.

There are various textbooks where L\'evy processes are discussed and where the above results may be either found directly or concluded from (e.g., \cite{a2009,B1996,ky2006,S2005} and pages 19-34 of \cite{p2003}). There are quite a few different proofs in the literature for the generalized PK formula, mostly via the application of the Wiener-Hopf factorization (e.g., \cite{H1977, Z1964}), weak convergence (e.g., \cite{CS2010,SW2002}) and martingales (e.g., \cite{kw1992}), but this is not the scope here.

For some recent work on the distribution (rather than Laplace-Stieltjes transform) of $\alpha$-stable L\'evy processes see \cite{BDP2008,K2011,Mi2011} and further references therein. For the case with phase type upward jumps (and general negative jumps) see \cite{Mo2002}. For some other results see also \cite{C2010}. These papers also include an extensive list of references to texts and further related literature.

One is not required to be a L\'evy process expert to read this paper and all the knowledge which is needed for what follows is covered above. In particular, this may easily be taught in any course where L\'evy processes are touched upon.

\section{The case of an independent sum of a Brownian motion and a subordinator}
Let $X$ be a L\'evy process which is an independent sum of a Brownian motion and a subordinator having a finite mean. In this case, as seen in the previous section, the (Laplace-Stieltjes) exponent can be written in the form
\begin{equation}
\phi(\alpha)=c\alpha+\frac{\sigma^2\alpha^2}{2}-\int_{(0,\infty)}(1-e^{-\alpha x})\nu(dx)\ ,
\end{equation}
where $\bar\nu\equiv\int_{(0,\infty)}x\nu(dx)=\int_0^\infty\nu(x,\infty)dx<c$. Denoting \begin{equation}
F_e(\alpha)=\int_0^\infty e^{-\alpha x}\frac{\nu(x,\infty)}{\bar\nu}dx\ ,
 \end{equation}
 $\rho=\frac{\bar\nu}{c}<1$ and $\lambda=\frac{2c}{\sigma^2}$, we observe that the exponent may be rewritten like this
\begin{equation}
\phi(\alpha)=c\alpha\left(1+\frac{\alpha}{\lambda}-\rho F_e(\alpha)\right)
\end{equation}
and in particular $\phi'(0)=c(1-\rho)$
so that the generalized PK formula has the form
\begin{equation}
\frac{\alpha\phi'(0)}{\phi(\alpha)}=\frac{1-\rho}{1+\frac{\alpha}{\lambda}-\rho F_e(\alpha)}
\end{equation}
If $\sigma^2=0$ ($\lambda=\infty$) and $\nu(0,\infty)<\infty$ then one obtains the well known (original) PK formula for the M/G/1 queue. When $\nu(0,\infty)=\infty$ it is interesting to observe that exactly the same formula is valid without change only that $F_e$ is then the stationary excess life distribution associated with the jumps of the subordinator and unlike in the renewal process setting, its density approaches $\infty$ in the neighborhood of zero. If the L\'evy measure is zero then $\rho=0$ and the formula becomes the Laplace-Stieltjes transform (LST) of an exponential distribution with rate $\lambda$. This is also well known to be the LST of the stationary distribution associated with a one dimensional reflected Brownian motion with negative drift. The interesting discovery is that these two results can be unified. To see this, note that
\begin{eqnarray}
\frac{\alpha\phi'(0)}{\phi(\alpha)}&=&\frac{\lambda}{\lambda+\alpha}
\frac{1-\rho}{1-\rho F_e(\alpha)\frac{\lambda}{\lambda+\alpha}}\n \\ \nonumber
&=&\sum_{n=0}^\infty (1-\rho)\rho^n F_e^n(\alpha)\left(\frac{\lambda}{\lambda+\alpha}\right)^{n+1}
\end{eqnarray}
and so we have the following result.
\begin{Theorem} \label{Th:main}
For a L\'evy process with no negative jumps satisfying $\int_{(0,1]}x\nu(dx)<\infty$, $\phi'(0)>0$ and with the notations defined above,
let $N\sim \mbox{G}(1-\rho)$ in the sense that $P[N=n]=(1-\rho)\rho^n$. Let $X_0,X_1,X_2,\ldots\sim \exp(\lambda)$ ($=0$ for $\lambda=\infty$) and $Y_1,Y_2,\ldots\sim F_e$ and assume that $N$, $X_0,X_1,\ldots$, $Y_1,Y_2,\ldots$ are all independent. Then $\frac{\alpha\phi'(0)}{\phi(\alpha)}$ is the LST of the following random variable
\begin{equation}
X_0+\sum_{n=1}^N(X_n+Y_n)
\end{equation}
where an empty sum is zero.
\end{Theorem}
When $\bar\nu\downarrow 0$ then we are left with $X_0$ as expected (the Brownian motion case) and when $\sigma^2\downarrow 0$ we are left with $\sum_{n=1}^NY_n$ also as expected (distributional form of the PK formula).

An interesting special case occurs when the jumps of the L\'evy process have a phase-type distribution. Then the residual life also has a phase-type distribution and thus also $X_i+Y_i$ and $X_0+\sum_{n=1}^N(X_n+Y_n)$. This a special case of the results reported in \cite{Mo2002} but with an easier proof but also a far less complicated setup.

It is very easy to check that also a converse holds.
\begin{Theorem}\label{Th:converse}
Assume that $0<p\le 1$, $0<\lambda\le\infty$, $f$ is nonnegative, nonincreasing with $\int_0^\infty f(y)dy=1$ (possibly with $f(x)\r\infty$ as $x\downarrow0$) and denote $F(x)=\int_0^xf(y)dy$. Let $N\sim\mbox{G}(p)$, $X_0,X_1,\ldots\sim \exp(\lambda)$ ($=0$ for $\lambda=\infty$) and $Y_1,Y_2,\ldots\sim F$ where all random variables are independent. Then there exists a L\'evy process which is an independent sum of a Brownian motion and a subordinator having a negative mean for which $\frac{\alpha\phi'(0)}{\phi(\alpha}$ is the LST of
\begin{equation}
X_0+\sum_{n=1}^N(X_n+Y_n)\ ,
\end{equation}
where an empty sum is zero. This L\'evy process is unique up to a constant time scale.
\end{Theorem}

\proof To see this we can simply perform reverse engineering. First assume without loss of generality that $f$ is right continuous, otherwise we take its right continuous version which gives the same $F$. Let $\nu((x,\infty))=\beta f(x)$ for any constant $\beta>0$ and thus $\nu((a,b])=\beta(f(a)-f(b)$, which uniquely characterizes $\nu$. Note that
\begin{equation}
\int_{(0,\infty)}x\nu(dx)=\int_0^\infty\nu((x,\infty))dx=\int_0^\infty \beta f(x)dx=\beta<\infty
\end{equation}
as required. Recalling that $\rho=\frac{\bar\nu}{c}$, that necessarily $p=1-\rho$ and since $\bar\nu=\beta$ we must set
$c=\frac{\beta}{1-p}$. Now, if $\lambda=\infty$ we set $\sigma^2=0$ and otherwise, from $\lambda=\frac{2c}{\sigma^2}$ we set $\sigma^2=\frac{2\beta}{\lambda(1-p)}$.
Finally,
\begin{equation}
\phi'(0)=c-\bar\nu=\frac{\beta}{1-p}-\beta>0
\end{equation}
so that all requirements are met.

The fact that this L\'evy process is unique up to a constant time scale is evident since if
\begin{equation}
\frac{\alpha\phi'(0)}{\phi(\alpha)}=\frac{\alpha\psi'(0)}{\psi(\alpha)}
\end{equation}
then necessarily $\psi(\alpha)=\gamma\phi(\alpha)$ for $\gamma=\frac{\psi'(0)}{\phi'(0)}$, so that if $X$ is a L\'evy process with exponent $\phi$ and $Y$ with $\psi$, then $\{X_{\gamma t}|\ t\ge0\}$ is distributed like $Y$.
\qed

\section{The general spectrally positive case}

The fact that $\frac{\alpha\phi'(0)}{\phi(\alpha)}$ is the LST of some proper distribution for any L\'evy process with no negative jumps and a negative mean is deduced indirectly from its various proofs which always exploit the connection with either the supremum or the reflected process. However, we are not aware of a more explicit and direct way of showing this than Theorems~\ref{Th:main},\ref{Th:converse} together with the following.

Observe that if $\int_{(0,1]}x\nu(dx)=\infty$ then the L\'evy process is not a sum of a Brownian motion and a subordinator. In this case (since $\phi'(0)>0$ and thus $\int_{(1,\infty)}x\nu(dx)<\infty$), we recall from (\ref{ea:mu}) in Section~\ref{Sec:intro} that
\begin{equation}
\phi(\alpha)=\mu\alpha+\frac{\sigma^2\alpha^2}{2}+\int_{(0,\infty)}\left(e^{-\alpha x}-1+\alpha x\right)\nu(dx)
\end{equation}
and that $\phi'(0)=\mu$. Now, with
\begin{equation}
\phi_\epsilon(\alpha)=\mu\alpha+\frac{\sigma^2\alpha^2}{2}+\int_{(\epsilon,\infty)}\left(e^{-\alpha x}-1+\alpha x\right)\nu(dx)
\end{equation}
for $\epsilon>0$, we clearly have that $\phi_\epsilon(\alpha)\r \phi(\alpha)$ as $\epsilon\downarrow 0$, that $\phi'_\epsilon(0)=\phi'(0)=\mu$ and thus
\begin{equation}
\lim_{\epsilon\downarrow0}\frac{\alpha\phi'_\epsilon(0)}{\phi_\epsilon(\alpha)}=\frac{\alpha \phi'(0)}{\phi(\alpha)}\ .
\end{equation}
and that $\frac{\alpha\phi'(0)}{\phi(\alpha)}\r1$ as $\alpha\downarrow 0$.

Now $\phi_\epsilon$ is the exponent of a L\'evy process which is an independent sum of a Brownian motion and a a compound Poisson process for which we have identified the distribution associated with $\frac{\alpha\phi'_\epsilon(0)}{\phi_\epsilon(\alpha)}$ in Theorem~\ref{Th:main}. From this the following is immediate.
\begin{Theorem}
 A distribution has the generalized PK LST for some L\'evy process with no negative jumps and a negative mean if and only if it belongs to the closure of the family of distributions defined in Theorem~\ref{Th:converse}.
\end{Theorem}


\begin{thebibliography}{99}
\bibitem{a2009} Applebaum, D. (2009). {\em L\'evy Processes and Stochastic Calculus, 2nd ed.}, Cambridge University Press.
\bibitem{BDP2008} Bernyk, V., Dalang, R. C. and Peskir G. (2008). The law of the supremum of a stable L\'evy process with no negative jumps. {\em Ann. Probab.} {\bf 36}, 1777-1789.
\bibitem{C2010} Chaumont, L. (2010). On the law of the supremum of L\'evy processes. arXiv:1011.4151v1 [math.PR].
\bibitem{B1996} Bertoin, J. (1996). {\em L\'evy Processes}, Cambridge University Press.
\bibitem{CS2010} Czysto{\l}owski, M. and W. Szczotka. (2010). Queueing approximation of suprema of spectrally positive L\'evy process. {\em Queueing Systems} {\bf 64}, 305-323.
    \bibitem{H1977} Harrison, J. M. (1977). The supremum distribution of a Levy process with no negative jumps. {\em Adv. Appl. Probab.} {\bf 9}, 417-422.
\bibitem{kw1992} Kella, O. and W. Whitt. (1992). Useful martingales for stochastic storage processes with L\'evy input. {\em J. Appl. Probab.} {\bf 29}, 396-403.
\bibitem{ky2006} Kyprianou, A. E. (2006). {\em Introductory Lectures on Fluctuations of L\'evy Processes with Applications}, Springer.
\bibitem{K2011} Kwa\'snicki, M., Ma{\l}ecki, J. and M. Ryznar. (2011). Suprema of L\'evy processes. arXiv:1103.0935v1 [math.PR].
\bibitem{Mi2011} Michna, Z. (2011). Formula for the supremum distribution of a spectrally positive $\alpha$-stable L\'evy process. {\em Stat. Probab. Letters} {\bf 81}, 231-235.
\bibitem{Mo2002} Mordecki, E. (2002). The distribution of the maximum of the L\'evy process with positive
jumps of phase-type. {\em Theory Stoch. Process.} {\bf 8}, 309-316.
\bibitem{p2003} Protter, P. E. (2003). {\em Stochastic Integration and Differential Equations, 2nd ed.}, Springer.
\bibitem{S2005} Sato, K. I. (1999). {\em L\'evy Processes and Infinitely Divisible Distributions}, Cambridge University Press.

\bibitem{SW2002} Szczotka, W. and W. A. Woyczy\'nski. (2003). Distributions of suprema of L´evy processes
via heavy traffic invariance principle. {\em Probab. Math. Stat.} {\bf 23}, 251-272.

\bibitem{Z1964} Zolotarev V.M. (1964) The first passage time to a level and the behavior at
infinity of processes with independent increments. {\em Theory Probab. Appl.} {\bf 9}, 653-
661.


\end{thebibliography}
\end{document}